\documentclass[a4paper]{amsart}
\usepackage{latexsym,amssymb,amsmath,amsfonts,amsthm}
\usepackage{epsf,psfrag,epsfig}
\usepackage[english]{babel}
\usepackage[latin1]{inputenc}

\newtheorem{thm}{Theorem}
\newtheorem{prop}[thm]{Proposition}

\newtheorem{cor}[thm]{Corollary}

\newcommand{\Ref}[1]{(\ref{#1})}

\newcommand{\parfrac}[2]{\left( \frac{#1}{#2} \right)}
\newcommand{\titre}[1]{\medskip  \noindent \textbf{#1}}
\newcommand{\ite}{\smallskip \noindent $\bullet$ }

\newcommand{\al}{\alpha}
\newcommand{\be}{\beta}
\newcommand{\si}{\sigma}
\newcommand{\la}{\lambda}

\newcommand{\ZZ}{\mathbb{Z}}

\newcommand{\pp}{\textbf{p}}

\newcommand{\mA}{\mathcal{A}}
\newcommand{\mB}{\mathcal{B}}
\newcommand{\mC}{\mathcal{C}}
\newcommand{\mD}{\mathcal{D}}
\newcommand{\mT}{\mathcal{T}}
\newcommand{\mU}{\mathcal{U}}
\newcommand{\mR}{\mathcal{R}}

\newcommand{\gA}{A}
\newcommand{\gB}{B}

\title{A Bijection between well-labelled positive paths and matchings}

\author[O. Bernardi]{Olivier Bernardi}
\address{O. Bernardi:  CNRS, Département de Mathématiques, Université Paris-Sud, 91405 Orsay, France}
\email{olivier.bernardi@math.u-psud.fr}

\author[B. Duplantier]{Bertrand Duplantier}
\address{B. Duplantier: Institut de Physique Théorique, Orme des Merisiers, CEA Saclay, 91191 Gif-sur-Yvette, France}
\email{bertrand.duplantier@cea.fr}

\author[P. Nadeau]{Philippe Nadeau}
\address{P. Nadeau: Fakult\"{a}t f\"{u}r Mathematik, Universit\"{a}t Wien, Nordbergstra{\ss}e 15, A-1090 Vienna, Austria}
\email{philippe.nadeau@univie.ac.at}

\date{\today}

\begin{document}
\maketitle

\begin{abstract}
A \emph{well-labelled positive path} of size $n$ is a pair $(\pp,\si)$ made of a word $\pp=p_1p_2\ldots p_{n-1}$ on the alphabet $\{-1, 0,+1\}$ such that $\sum_{i=1}^j p_i \geq 0$ for all $j=1\ldots n\!-\!1$, together with a permutation $\si=\si_1\si_2\ldots\si_n$ of $\{1,\ldots,n\}$ such that $p_i=-1$ implies $\si_{i}<\si_{i+1}$, while $p_i=1$ implies $\si_{i}>\si_{i+1}$.
We establish a bijection between well-labelled positive paths of size $n$ and matchings (i.e. fixed-point free involutions) on $\{1,\ldots,2n\}$. This proves that the number of well-labelled positive paths is $(2n-1)!!\equiv(2n-1)\cdot(2n-3)\cdots3\cdot 1$. By specialising our bijection, we also prove that the number of permutations of size $n$ such that each prefix has no more ascents than descents is $[(n-1)!!]^2$ if $n$ is even and $n!!\,(n-2)!!$ otherwise. 

It is shown in \cite{BDN:Freely-jointed-chain} that well-labelled positive paths of size $n$ are in bijection with a collection of $n$-dimensional subpolytopes partitioning the polytope $\Pi_n$ consisting of all points $(x_1,\ldots,x_n)\in [-1,1]^n$ such that $\sum_{i=1}^j x_i\geq 0$ for all $j=1 \ldots n$. Given that the volume of each subpolytope is $1/n!$, our results prove combinatorially that the volume of $\Pi_n$ is $\frac{(2n-1)!!}{n!}$.
\end{abstract}

\section{Introduction}
\begin{figure}[ht!]\begin{center}
\input{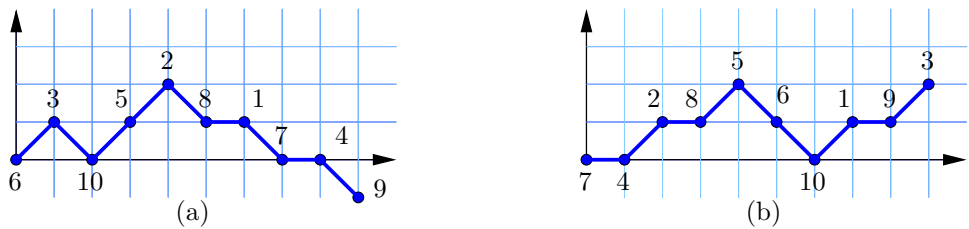}\caption{(a) A well-labelled Motzkin path. (b) A well-labelled positive path.}\label{fig:exp-labelled-path} \end{center}\end{figure}

A \emph{well-labelled path} of size $n$ is a pair $(\pp,\si)$ made of a word $\pp=p_1p_2\ldots p_{n-1}$ on the alphabet $\{-1, 0,+1\}$, together with a permutation $\si=\si_1\si_2\ldots\si_n$ of $[n]\equiv\{1,\ldots,n\}$ such that $p_i=-1$ implies $\si_{i}<\si_{i+1}$, while $p_i=1$ implies $\si_{i}>\si_{i+1}$.
We shall represent a path $(\pp,\si)$ as a lattice path on $\ZZ^2$ starting at (0,0) and made of steps $(1,p_i)$ for $i=1\ldots n\!-\!1$ together with the label $\si_i$ on the $i$th lattice point of the path for $i=1\ldots n$. For instance, two well-labelled paths of size 10 are represented in Figure~\ref{fig:exp-labelled-path}. A well-labelled path $(\pp,\si)$ is said \emph{Motzkin} (resp. \emph{positive}) if $\displaystyle \sum_{i=1}^j p_i \geq 0$ for all $j=1\ldots n\!-\!2$ and $\displaystyle\sum_{i=1}^{n-1} p_i = -1$ (resp. $\displaystyle \sum_{i=1}^j p_i \geq 0$ for all $j=1 \ldots n\!-\!1$). \\

In this paper, we define a bijection between well-labelled Motzkin paths of size $n+1$ and \emph{matchings} (fixed-point free involutions) on $[2n]$. We then define a closely related bijection between well-labelled positive paths of size $n$, and matchings on $[2n]$.  This proves that these sets of paths are counted by $(2n-1)!!\equiv(2n-1)\cdot(2n-3)\cdots3\cdot 1$. Both bijections also allow for a refined enumeration of well-labelled paths $(\pp,\si)$ according to the number of \emph{horizontal steps} (the number of letters 0 in $\pp$). More precisely, we  show that the number of well-labelled positive paths of size $n$ with $k$ horizontal steps is $\binom{n}{k}\binom{n-1}{k}\, k!\,[(n-k-1)!!]^2$ if $n-k$ is even and $\binom{n}{k}\binom{n-1}{k}\,k!\,(n-k)!!\,(n-k-2)!!$ otherwise. Similarly, the number of well-labelled Motzkin paths of size $n$ with $k$ horizontal steps is $ \binom{n}{k}\binom{n-2}{k}\,k!(n-k-1)!!(n-k-3)!!$ if $n-k$ is even and 0 otherwise.
Observe that well-labelled paths $(\pp,\si)$ without horizontal steps are completely determined by the permutation $\si$. Indeed, in this case the word $\pp$ encodes the \emph{up-down sequence} of the permutation $\si$.
Hence, by specialising our results to paths with no horizontal steps (i.e. $k=0$), we enumerate permutations whose up-down sequence belong to a certain family. For instance, we prove that the number of permutations of size $n$ such that each prefix has no more ascents than descents is $[(n-1)!!]^2$ if $n$ is even and $n!!\,(n-2)!!$ otherwise.  We also prove that 
 the number of permutations of size $n$ having one more ascent than descent but such that each prefix has no more ascents than descents is $(n-1)!!\, (n-3)!!$ if $n$ is even and 0 otherwise.
These enumerative results contrast with those in \cite{Carlitz:prescribed-patterns,deBruijn:perm-given-updown,Foulkes:enum-perm-updown,Niven:comb-problem} by the fact that we consider here a family of admissible up-down sequences rather than a single sequence.\\

Well-labelled positive paths appeared recently in a problem concerning the evaluation of the volume of the $n$-dimensional polytope $\Pi_n$ made of the points $(x_1,\ldots,x_n)$ in $[-1,1]^n$ such that $\sum_{i=1}^j x_i\geq 0$ for all $j=1 \ldots n$.  Indeed, it was shown in \cite{BDN:Freely-jointed-chain} that the set of well-labelled positive paths of size $n$ is in bijection with a set of $n$-dimensional subpolytopes forming a partition of $\Pi_n$ and this was our original motivation for studying well-labelled paths.  Given that the volume of each subpolytope is $1/n!$, our results prove combinatorially that the volume of $\Pi_n$ is $\frac{(2n-1)!!}{n!}$.\\

The paper is organised as follows. In Section~\ref{section:decomposition}, we define a recursive decomposition of well-labelled positive and Motzkin paths. We then translate these decompositions in terms of generating functions. For Motzkin paths, solving the generating function equation shows that the number of well-labelled Motzkin paths of size $n+1$ is $(2n-1)!!$. From this, a simple induction shows that the number of well-labelled positive paths of size $n$ is also $(2n-1)!!$. In Section~\ref{section:bijection}, we use the recursive decomposition of paths in order to define bijections between well-labelled positive paths, well-labelled Motzkin paths and matchings. One step of these bijections uses a construction of Chen \cite{Chen:bij-algo-trees} between labelled binary trees and matchings. Lastly in Section~\ref{section:counting}, we use our bijections to count well-labelled positive and Motzkin paths according to their number of horizontal steps. Specialising this results to the paths with no horizontal steps, we enumerate permutations whose up-down sequence belong to certain families mentioned above.\\


\section{Decomposition of well-labelled paths} \label{section:decomposition}
In this section, we define a recursive decomposition of the class $\mA$ of well-labelled Motzkin paths and the class $\mB$ of well-labelled positive paths. We then translate these equations in terms of generating functions and obtain our first counting results. \\

We denote respectively by $\mA_n$ and  $\mB_n$  the sets of paths of size $n$ in $\mA$ and $\mB$. We denote respectively by $a_n$ and $b_n$ the cardinality of $\mA_n$ and $\mB_n$ and by
$$\gA(z)=\sum_{n\geq 0}\frac{a_n}{n!}z^n ~~\textrm{ and }~~ \gB(z)=\sum_{n\geq 0}\frac{b_n}{n!}z^n$$
the corresponding exponential generating functions. Observe that $a_0=a_1=0$ and $b_0=0$. The following notation will be useful for \emph{relabelling} objects: given a set  $I$ of $n$ integers, we denote by $\la_I$ the order preserving bijection from $[n]$ to $I$ (and by $\la^{-1}_I$ the inverse bijection).\\

\subsection{Decomposition of well-labelled Motzkin paths.}
\begin{figure}[ht!]\begin{center}
\input{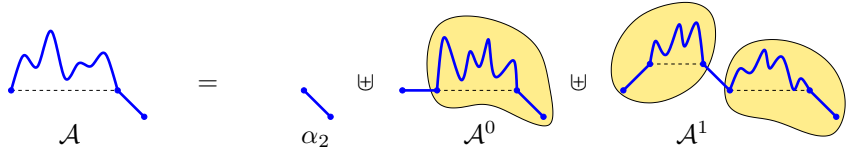}\caption{Recursive decomposition of well-labelled Motzkin paths.}\label{fig:decomposition-Motzkin} \end{center}\end{figure}

We first define a recursive decomposition of the class $\mA$ of well-labelled Motzkin paths. For $i\in\{-1,0,1\}$, we denote by $\mA^i$ (resp. $\mA^i_n$) the set of paths $(\pp,\si)$ in $\mA$ (resp. $\mA_n$) such that $p_1=i$. Observe that $\mA^{-1}$ is made of a single element $\al_2$ of size 2. The decomposition
$$\mA=\{\al_2\}\uplus\mA^{0}\uplus\mA^{1}$$
is illustrated by Figure~\ref{fig:decomposition-Motzkin} and the following proposition reveals its recursive nature.

\begin{prop}\label{prop:A} For any positive integer $n$,
\begin{itemize}
\item the set $\mA_{n}^0$ is in bijection with the set $[n]\times \mA_{n-1}$,
\item the set $\mA_{n}^1$ is in bijection with the set $\mC_n$ made of all unordered pairs $\{(I',P'),(I'',P'')\}$ such that $I'\subseteq [n]$, $I''=[n]\setminus I'$ and $P'$, $P''$ are well-labelled Motzkin paths of respective size $|I'|$ and $|I''|$.
\end{itemize}
\end{prop}

\begin{proof}
\ite For any path $(\pp,\si)$ in $\mA_{n}^0$ one obtains a pair $(k,(\pp',\si'))$ in $[n]\times \mA_{n-1}$ by setting $k=\si_1$, $\pp'=p_2\ldots p_{n-1}$ and $\si'=\si'_1\ldots \si'_{n-1}$ where $\si'_i=\la^{-1}_{[n]\setminus\{k\}}(\si_{i+1})$ for $i=1\ldots n\!-\!1$. This is clearly a bijection between $\mA_{n}^0$ and $[n]\times \mA_{n-1}$.\\
\ite Observe that the set $\mC_n$ is trivially in bijection with the set $\mC_n'$ of \emph{ordered} pairs $((I',P'),(I'',P''))$ such that $\la_{I'}(\si'_1)<\la_{I''}(\si''_1)$. Let $(\pp,\si)$ be a path in $\mA_n^1$ and let $k$ be the least integer such that $\sum_{i=1}^k p_i=0$ (observe that $p_k=-1$). We define an element $\phi(\pp,\si)=((I',(\pp',\si')),(I'',(\pp'',\si'')))$ of $\mC_n'$ by setting
\begin{itemize}
\item $I'=\{\si_1,\ldots,\si_k\}$ and $I''=\{\si_{k+1},\ldots,\si_n\}$,
\item $\pp'=p'_1\ldots p'_{k-1}$  and $\pp''=p''_1\ldots p''_{n-k-1}$, where $p'_i=-p_{k-i}$ and  $p_i''=p_{k+i}$,
\item $\si'=\si'_1\ldots \si'_k$ and $\si''=\si''_1\ldots \si''_{n-k}$, where $\si'_i=\la^{-1}_{I'}(\si_{k+1-i})$ and  $\si''_i=\la^{-1}_{I''}(\si_{k+i})$.
\end{itemize}
The mapping $\phi$ is clearly a bijection between the sets $\mA_{n}^1$ and $\mC_n'$, which concludes the proof.
\end{proof}

\begin{cor}\label{cor:A}
The generating function $\gA(z)$ of well-labelled Motzkin paths satisfies
\begin{eqnarray}\label{eq:A}
\gA(z)~=~\frac{z^2}{2}+z\gA(z)+\frac{\gA(z)^2}{2}.
\end{eqnarray}
\end{cor}

\begin{proof} For $i\in\{-1,0,1\}$, we denote by $a_n^i$ the cardinality of $\mA_n^i$ and by $\gA^i(z)=\sum_{n\geq 0}\frac{a^i_n}{n!}z^n$ the corresponding generating function. The partition $\mA=\{\al_2\}\uplus\mA^{0}\uplus\mA^{1}$ gives
$$\gA(z)=\frac{z^2}{2}+\gA^0(z)+\gA^1(z).$$
Moreover, the bijection between $\mA^{0}_n$ and $[n]\times \mA_{n-1}$ gives $a^0_{n}=n\, a_{n-1}$, hence $\gA^0(z)=z\gA(z)$ while the correspondence between $\mA^{1}_n$ and $\mC_n$ gives $\displaystyle a^1_{n}=\frac{1}{2}\sum_{k=0}^n \binom{n }{k}a_ka_{n-k}$, hence $\gA^1(z)=\frac{\gA(z)^2}{2}.$
\end{proof}

By solving Equation \Ref{eq:A} (and using the fact that $a_0=0$), one gets 
\begin{eqnarray}\label{eq:Abis}
A(z)=1-z-\sqrt{1-2z}.
\end{eqnarray}
One can extract the coefficient $a_n$ either directly from this expression of $A(z)$ or by applying Lagrange inversion formula to the series $C(z)=\gA(z)/z$. Indeed, Equation~\Ref{eq:A} gives $C(z)=z\frac{(1+C(z))^2}{2}$, hence
$$a_{n+1}=(n+1)!\displaystyle[z^{n}]C(z)=\frac{(n+1)!}{n}[x^{n-1}]\parfrac{(1+x)^{2}}{2}^n=\frac{(2n)!}{2^nn!}=(2n-1)!!.$$

We will now determine the number $b_n$ of well-labelled positive paths of size $n$. This can be done by exploiting a bijection  between $\mB_n\times[n+1]\times\{0,1\}$ and $\mB_{n+1}\uplus \mA_{n+1}$ obtained by \emph{adding one step} to a positive path. The bijection is as follows: given a well-labelled positive path $(\pp,\si)$ of size $n$, an integer $k$ in $[n+1]$ and an integer $b$ in $\{0,1\}$, we define the labelled path $(\pp',\si')=\psi((\pp,\si),k,b)$ by setting
\begin{itemize}
\item $\si'=\si'_1\ldots\si'_{n+1}$, where   $\si_i'=\la_{[n+1]\setminus\{k\}}(\si_i)$  for $i=1\ldots n$ and $\si'_{n+1}=k$,
\item $\pp'=p'_1\ldots p'_n$, where $p_i'=p_i$ for $i=1\ldots n\!-\!1$ and  $p'_n$ is equal to $b-1$ if $\si'_{n+1}>\si'_n$ and equal to $b$ otherwise.
\end{itemize}
Observe that the path $(\pp',\si')$ is well-labelled (by the choice of the step $p'_n$) and is either positive or Motzkin (since $(\pp,\si)$ is positive). Moreover, the mapping $\psi$ is a bijection between $\mB_n\times[n+1]\times\{0,1\}$ and $\mB_{n+1}\uplus \mA_{n+1}$ showing that
\begin{eqnarray}\label{eq:anbn}
2(n+1)b_n=b_{n+1}+a_{n+1}\text{~~for all~}n\geq 0.
\end{eqnarray}
Since $a_{n+1}=(2n-1)!!$ a simple induction shows that $b_{n}=(2n-1)!!$ and proves the following.

\begin{prop}\label{prop:double-factorial}
The number $a_{n+1}$ of well-labelled Motzkin paths of size $n+1$ and the number $b_n$ of well-labelled positive paths of size $n$ are both equal to $(2n-1)!!$.\\
\end{prop}

\subsection{Decomposition of well-labelled positive paths.}
\begin{figure}[ht!]\begin{center}
\input{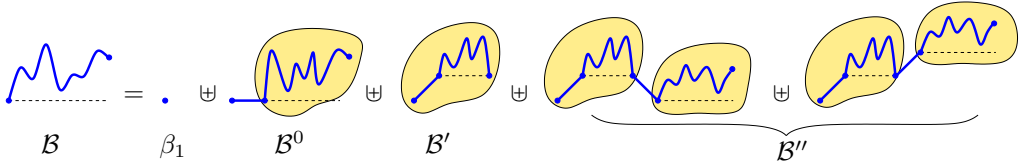}\caption{Recursive decomposition of well-labelled positive paths.}\label{fig:decomposition-PrefMotzkin} \end{center}\end{figure}

We now define a recursive decomposition of the class $\mB$ of well-labelled positive paths. We denote by $\be_1$ the well-labelled path of size $1$ and for $i\in\{0,1\}$, we denote by $\mB^i$ the set of paths $(\pp,\si)$ in $\mB$ of size at least 2 satisfying $p_1=i$. For a path $(\pp,\si)$ of size $n$ in $\mB^1$, we consider the greatest integer $k\leq n$ such that 
$\sum_{i=1}^{j-1} p_i\geq 1$ for all $j=2\ldots k\!-\!1$ and $\sum_{i=1}^{k-1} p_i=1$. We denote by $\mB'$ the subset of paths in $\mB^1$ such that $k=n$ and we denote $\mB''=\mB^1\setminus \mB'$ the complement. We also denote by $\mB^0_n$, $\mB'_n$ and $\mB''_n$ respectively the paths of size $n$ in  $\mB^0$, $\mB'$ and $\mB''$.
The partition
$$\mB=\{\be_1\}\uplus\mB^0\uplus\mB'\uplus\mB''$$
is illustrated by Figure~\ref{fig:decomposition-PrefMotzkin} and the following proposition reveals its recursive nature.

\begin{prop}
\label{prop:B}
For any positive integer $n$, 
\begin{itemize}
\item the set $\mB^0_n$ is in bijection with the set $[n]\times \mB_{n-1}$,
\item the set $\mB'_n$ is in bijection with the class $\mA_n$ of well-labelled Motzkin paths,
\item the set $\mB''_n$ is in bijection with the set $\mD_n$ made of the ordered pairs $((I',P'),(I'',P''))$ such that $I'\subseteq [n]$, $I''=[n]\setminus I'$, $P'$ is a well-labelled Motzkin path of size $|I'|$ and $P''$ is a well-labelled positive path of size $|I''|$.
\end{itemize}
\end{prop}

\begin{proof}
\ite For any path $(\pp,\si)$ in $\mB_{n}^0$ one obtains a pair $(k,(\pp',\si'))$ in $[n]\times \mB_{n-1}$ by setting $k=\si_1$, $\pp'=p_2\ldots p_{n-1}$ and $\si'=\si'_1\ldots \si'_{n-1}$ where $\si'_i=\la^{-1}_{[n]\setminus\{k\}}(\si_{i+1})$ for $i=1\ldots n\!-\!1$. This is clearly a bijection between $\mB_{n}^0$ and $[n]\times \mB_{n-1}$.\\
\ite A bijection between the sets  $\mB'_n$ and $\mA_n$ is obtained by reading the positive path backward: given a path $(\pp,\si)$ in $\mB_{n}'$ one obtain a path $(\pp',\si')$ in $\mA_n$ by setting $\si'=\si_n\ldots \si_1$ and $\pp'=p'_1\ldots p'_{n-1}$ where $p_i'=-p_{n-i}$ for $i=1\ldots n\!-\!1$. This is clearly a bijection.\\
\ite Let $(\pp,\si)$ be a path in $\mB''_n$ and let $k< n$ be the greatest integer such that $\sum_{i=1}^{j-1} p_i\geq 1$ for all $j=2\ldots k\!-\!1$ and $\sum_{i=1}^{k-1} p_i=1$. We define a pair $((I',(\pp',\si')),(I'',(\pp'',\si'')))=\phi(\pp,\si)$ by setting
\begin{itemize}
\item $I'=\{\si_1,\ldots,\si_k\}$ and $I''=\{\si_{k+1},\ldots,\si_n\}$,
\item $\pp'=p'_1\ldots p'_{k-1}$  and $\pp''=p''_1\ldots p''_{n-k-1}$, where $p'_i=-p_{k-i}$ and  $p_i''=p_{k+i}$,
\item $\si'=\si'_1\ldots \si'_k$ and $\si''=\si''_1\ldots \si''_{n-k}$, where $\si'_i=\la^{-1}_{I'}(\si_{k+1-i})$ and  $\si''_i=\la^{-1}_{I''}(\si_{k+i})$.
\end{itemize}
We first want to prove that $((I',(\pp',\si')),(I'',(\pp'',\si'')))$ is in $\mD_n$. It is clear that $(\pp',\si')$ and $(\pp'',\si'')$ are well-labelled paths and moreover, $(\pp',\si')$ is a Motzkin path. It remains to prove that $(\pp'',\si'')$ is a positive path. Observe that the step $p_k$ is non-zero otherwise it contradicts the maximality of $k$. If $p_k=-1$, then $(\pp'',\si'')$ is clearly positive because $(\pp,\si)$ is positive; and if  $p_k=+1$, then $(\pp'',\si'')$ is positive otherwise it would contradict the maximality of $k$.  Hence, $\phi$ is a mapping from $\mA_n$ to $\mD_n$. The bijectivity of $\phi$ is easy to check after observing that the step $p_k$ can be recovered: it is equal to $1$ if $\la_{I'}(\si'_1)<\la_{I''}(\si''_1)$ and to $-1$ otherwise.
\end{proof}

Proposition~\ref{prop:B} will allow to define a bijection between positive paths and matchings in the next section. It also leads to the following relation between the generating functions $A(z)$ and $B(z)$: 
\begin{equation}
 \gB(z)=z+z\gB(z)+\gA(z)+\gA(z)\gB(z),
 \label{eq:BA}
\end{equation}
which, by \Ref{eq:Abis}, gives $\displaystyle B(z)=\frac{1}{\sqrt{1-2z}}-1$. This result could also have been derived from the observation that $b_n=a_{n+1}$ implies $B(z)=A'(z)$.




\section{Bijections with matchings} \label{section:bijection}
Proposition~\ref{prop:double-factorial} suggests that the classes of paths $\mA_{n+1}$ and $\mB_n$ are both in bijection with matchings. The goal of this section is to describe such bijections. For this, we will introduce intermediate objects called \emph{labelled binary trees}.

\subsection{Bijections between well-labelled paths and labelled binary trees}
\begin{figure}[ht!]\begin{center}
\input{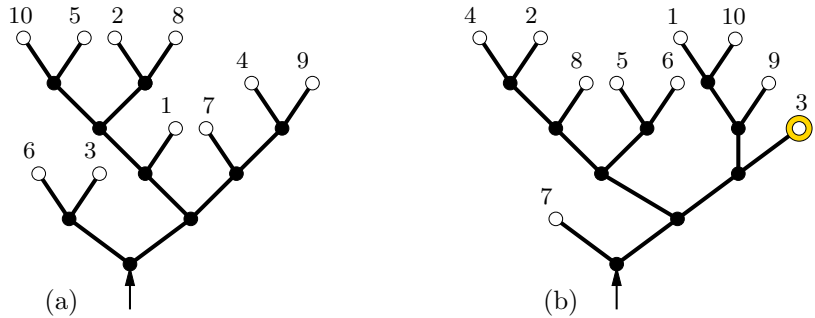}\caption{(a) A labelled binary tree. (b) A marked labelled binary tree}\label{fig:exp-labelled-trees} \end{center}\end{figure}

A \emph{labelled binary tree} of size $n$ is a rooted tree with $n$ leaves having $n$ different labels in $[n]$ and such that each (unlabelled) internal vertex has exactly two unordered children.  We call \emph{marked labelled binary tree} a labelled binary tree in which one of the (internal or external) vertices is marked. A binary tree and a marked binary tree are represented in Figure~\ref{fig:exp-labelled-trees}. We denote by  $\mT$ the set of labelled binary trees \emph{of size at least} 2 and we denote by $\mR$ the set of marked labelled binary trees. We will now show that the recursive descriptions of the classes $\mT$ and $\mR$ parallel those of the classes $\mA$ and $\mB$ and obtain bijections between $\mT$ and $\mA$ and between $\mR$ and $\mB$.\\

We use the following notation for \emph{relabelling} trees: if $\la$ is a bijection between two sets of integers $I$, $J$ and $\tau$ is a binary tree whose leaves have labels in $I$, then $\la(\tau)$ denotes the tree obtained from $\tau$ by replacing each leaf labelled $i\in I$, $i=1\ldots n$ by a leaf labelled $\la(i)\in J$.

\titre{Bijection between Motzkin paths and labelled trees.}
We denote by $\tau_2$ the unique labelled binary tree of size 2 and by $\mT^0$ (resp. $\mT^1$) the set of labelled binary trees made of the root-vertex and two subtrees, one of which (resp. none of which) is a leaf. For all integers $n\geq 2$, we denote respectively by $\mT_n$, $\mT_n^{0}$ and $\mT_n^{1}$ the set of trees of size $n$ in  $\mT$, $\mT^{0}$ and $\mT^{1}$. We now explicit the recursive nature of the partition
$$\mT=\{\tau_2\}\uplus \mT^0\uplus \mT^1.$$

\begin{prop}\label{prop:T}
For any integer $n> 2$,
\begin{itemize}
\item the set $\mT_{n}^0$ is in bijection with the set $[n]\times \mT_{n-1}$.
\item the set $\mT_{n}^1$ is in bijection with the set $\mU_n$ of unordered pairs $\{(I',\tau'),(I'',\tau'')\}$ such that $I'\subseteq [n]$, $I''=[n]\setminus I'$ and $\tau'$, $\tau''$ are labelled binary trees in $\mT$ of respective size $|I'|$ and~$|I''|$.
\end{itemize}
\end{prop}

\begin{proof} Let $n>2$.\\ 
\ite Given an integer $k\in [n]$ and a tree $\tau$ in $\mT_{n-1}$, one defines a tree $\tau'$ in $\mT_{n}^0$ as the tree made of a root-vertex, one subtree which is a leaf labelled $k$ and one subtree equal to $\la_{[n]\setminus\{k\}}(\tau)$ (which is not a leaf). This is a bijection between $[n]\times \mT_{n-1}$ and $\mT_{n}^0$.\\
\ite Given a pair $\{(I',\tau'),(I'',\tau'')\}$ in $\mU_n$, one defines a tree $\tau$ in $\mT_{n}^1$ as the tree made of a root-vertex, a subtree equal to $\la_{I'}(\tau')$ and a subtree equal to $\la_{I''}(\tau'')$. This is a bijection between $\mU_n$ and $\mT_{n}^1$.
\end{proof}

\textbf{Definition of bijection $\Phi$.} Comparing Propositions~\ref{prop:A} and~\ref{prop:T}, it is clear that one can define a recursive bijection $\Phi$ between the class $\mA$ of well-labelled Motzkin paths and the class $\mT$ of labelled binary trees.
We now summarise the recursive construction of the image of a well-labelled Motzkin path $(\pp,\si)$ by the bijection~$\Phi$:
\begin{itemize}
\item[(i)] If $(\pp,\si)=\al_2$, then $\Phi(\pp,\si)=\tau_2$.
\item[(ii)]  If $(\pp,\si)$ has size $n>2$ and $p_1=0$, then we set $\pp'=p_2\ldots p_{n-1}$ and $\si'=\si'_1\ldots \si'_{n-1}$ where $\si_i'=\la^{-1}_{[n]\setminus\{\si_1\}}(\si_{i+1})$ for $i=1\ldots k\!-\!1$. With these notations, we define  $\Phi(\pp,\si)$ as the tree made of a root-vertex, the subtree made of a leaf labelled $\si_1$ and the subtree $\la_{[n]\setminus\{\si_1\}}(\Phi(\pp',\si'))$.
\item[(iii)]  If $(\pp,\si)$ has size $n>2$ and $p_1=1$, then we consider the least integer $k$ such that $\sum_{i=1}^kp_i=0$ and we set (as in the proof of Proposition \ref{prop:A}):\\
\ite  $I'=\{\si_1,\ldots,\si_k\}$ and $I''=\{\si_{k+1},\ldots,\si_n\}$,\\
\ite $\pp'=p'_1\ldots p'_{k-1}$ and $\pp''=p''_1\ldots p''_{n-k-1}$, where $p'_i=-p_{k-i}$ and $p''_i=p_{k+i}$,\\
\ite $\si'=\si'_1\ldots \si'_k$ and $\si''=\si''_1\ldots \si''_{n-k}$, where $\si'_i=\la^{-1}_{I'}(\si_{k+1-i})$ and  $\si''_i=\la^{-1}_{I''}(\si_{k+i})$.\\
With these notations, we define $\Phi(\pp,\si)$ as the tree made of a root-vertex, the subtree $\la_{I'}(\Phi(\pp',\si'))$ and  the subtree $\la_{I''}(\Phi(\pp'',\si''))$.
\end{itemize}

For instance, the image of the Motzkin path represented in Figure~\ref{fig:exp-labelled-path}(a) by the mapping $\Phi$ is represented in Figure~\ref{fig:exp-labelled-trees}(a). From the definition of $\Phi$ and Propositions \ref{prop:A} and \ref{prop:T}, we have the following theorem:

\begin{thm} \label{thm:bij1}
For any positive integer $n>1$, the mapping $\Phi$ is a bijection between well-labelled Motzkin paths of size $n$ and labelled binary trees with $n$ leaves.
\end{thm}

\titre{Bijection between positive paths and marked trees.}
We now define a bijection $\Phi'$ between well-labelled positive paths and marked labelled binary trees. Before defining the bijection $\Phi'$, let us explain briefly what led us to consider \emph{marked} labelled binary trees. As seen in Section 2, the recursive decomposition of positive paths leads to consider
blocks corresponding to either positive paths or Motzkin paths.
The recursive relation is captured by Equation (\ref{eq:BA}) which can be written
$$B(z)=\tilde{A}(z)+\tilde{A}(z)B(z),$$
where $\tilde{A}(z)=z+A(z)$ is the series of unmarked binary trees
(of size $n\geq 1$), or equivalently, binary trees marked at their
root-vertex. This relation suggests that one can interpret $B(z)$ as
counting marked binary trees.\\

We denote by $\rho_1$ the marked labelled tree of size 1. We denote by $\mR'\subset \mR$ the set of marked trees of size at least 2 such that the marked vertex is the root. Clearly, this set is in bijection with the set $\mT$ of unmarked trees. We denote by $\mR^0$ (resp. $\mR''$) the set of marked trees made of a non-marked root and two subtrees, one of which (resp. none of which) is a non-marked leaf. For all integer $n>1$, we denote respectively by $\mR_n$, $\mR_n^{0}$, $\mR_n'$ and $\mR_n''$ the set of marked trees of size $n$ in   $\mR$, $\mR^{0}$, $\mR'$ and $\mR''$. We now explicit the recursive nature of the partition
$$\mR=\{\rho_1\}\uplus \mR^{0}\uplus\mR'\uplus\mR''.$$

\begin{prop}\label{prop:R} For all integer $n>1$,
\begin{itemize}
\item the set $\mR_n^0$ is in bijection with the set $[n]\times \mR_{n-1}$,
\item the set $\mR_n'$ is in bijection with $\mT_n$,
\item the set $\mR_n''$ is in bijection with the set $\mU_n$ of ordered pairs $((I',\tau'),(I'',\tau''))$ such that $I'\subseteq [n]$, $I''=[n]\setminus I'$, $\tau'$ is a non-marked labelled binary tree of size $|I'|$ and $\tau''$ is a marked labelled binary tree of size $|I''|$.
\end{itemize}
\end{prop}

The proof of Proposition~\ref{prop:R} is similar to the proof of Proposition~\ref{prop:T} and is omitted.\\


\textbf{Definition of bijection $\Phi'$.} Comparing Propositions~\ref{prop:B} and~\ref{prop:R}, it is clear that one can define a recursive bijection $\Phi'$ between the class $\mB$ of well-labelled paths and the class $\mR$ of marked labelled binary trees.
We now summarise the recursive construction of the image of a well-labelled positive path $(\pp,\si)$ by the bijection~$\Phi'$:
\begin{itemize}
\item[(i)] If $(\pp,\si)=\beta_1$ then $\Phi'(\pp,\si)$ is the marked tree $\rho_1$.
\item[(ii)] If $(\pp,\si)$ has size $n>1$ and $p_1=0$, then we define  $\pp'=p_2\ldots p_{n-1}$ and $\si'=\si'_1\ldots\si'_{n-1}$, where $\si'_i=\la^{-1}_{[n]\setminus\{\si_1\}}(\si_{i+1})$ for $i=1\ldots n\!-\!1$.  With these notations, we define $\Phi'(\pp,\si)$ as the tree made of a non-marked root-vertex, the subtree made of a non-marked leaf labelled $\si_1$ and the marked subtree $\la_{[n]\setminus\{\si_1\}}(\Phi'(\pp',\si'))$.
\item[(iii)] If $(\pp,\si)$ has size $n>1$ and $p_1=1$, then we consider the greatest integer $k\leq n$ such that
$\sum_{i=1}^{j-1} p_i\geq 1$ for all $j=2\ldots k\!-\!1$ and $\sum_{i=1}^{k-1} p_i=1$, and we set (as in the proof of Proposition \ref{prop:B}):\\
\ite $I'=\{\si_1,\ldots,\si_k\}$ and $I''=\{\si_{k+1},\ldots,\si_n\}$,\\
\ite $\pp'=p'_1\ldots p'_{k-1}$  and $\pp''=p''_1\ldots p''_{n-k-1}$, where $p'_i=-p_{k-i}$ and  $p_i''=p_{k+i}$,\\
\ite $\si'=\si'_1\ldots \si'_k$ and $\si''=\si''_1\ldots \si''_{n-k}$, where $\si'_i=\la^{-1}_{I'}(\si_{k+1-i})$ and  $\si''_i=\la^{-1}_{I''}(\si_{k+i})$.\\
If $k=n$ (that is, $(\pp,\si)$ is in $\mB'$), we define  $\Phi'(\pp,\si)$ as the marked tree obtained by marking the root-vertex of the unmarked tree $\Phi(\pp',\si')$ (note that $I''$, $\pp''$ and $\si''$ are empty in this case). Otherwise (that is, if $k<n$), we define $\Phi'(\pp,\si)$ as the marked tree made of a non-marked root-vertex, the non-marked subtree $\la_{I'}(\Phi(\pp',\si'))$ and the marked subtree $\la_{I''}(\Phi'(\pp'',\si''))$.
\end{itemize}

For instance, the image of the positive path represented in Figure~\ref{fig:exp-labelled-path}(b) by the mapping $\Phi'$ is represented in Figure~\ref{fig:exp-labelled-trees}(b).  From the definition of $\Phi'$ and Propositions \ref{prop:B} and \ref{prop:R}, we have the following theorem:

\begin{thm}
 The function $\Phi'$ is a bijection between well-labelled positive paths of size $n$ and marked labelled binary trees with $n$ leaves.
 \end{thm}

\noindent \textbf{Remark.} The \emph{final height} of a positive path $(\pp,\si)$ of size $n$ is $\sum_{i=1}^{n-1}p_i$. It is not hard to prove inductively that a positive path $(\pp,\si)$ has an even final height if and only if the mark of the image tree, $\rho=\Phi'(\pp,\si)$, is on a leaf. Indeed, if $(\pp,\si)=\beta_1$ the final height is 0, and the mark is on a leaf of the tree $\rho=\rho_1$; while if $(\pp,\si)$ is in $\mB'$ the final height is 1, and the mark is on an internal vertex (the root-vertex) of $\rho$. In the other cases ($(\pp,\si)\in\mB_0\uplus\mB''$), the mark is in a subtree $\rho'$ of $\rho$ corresponding to a path $(\pp',\si')$ having a final height of the same parity as $(\pp,\si)$.\\

\subsection{Bijections between labelled binary trees and matchings}
\begin{figure}[ht!]\begin{center}
\input{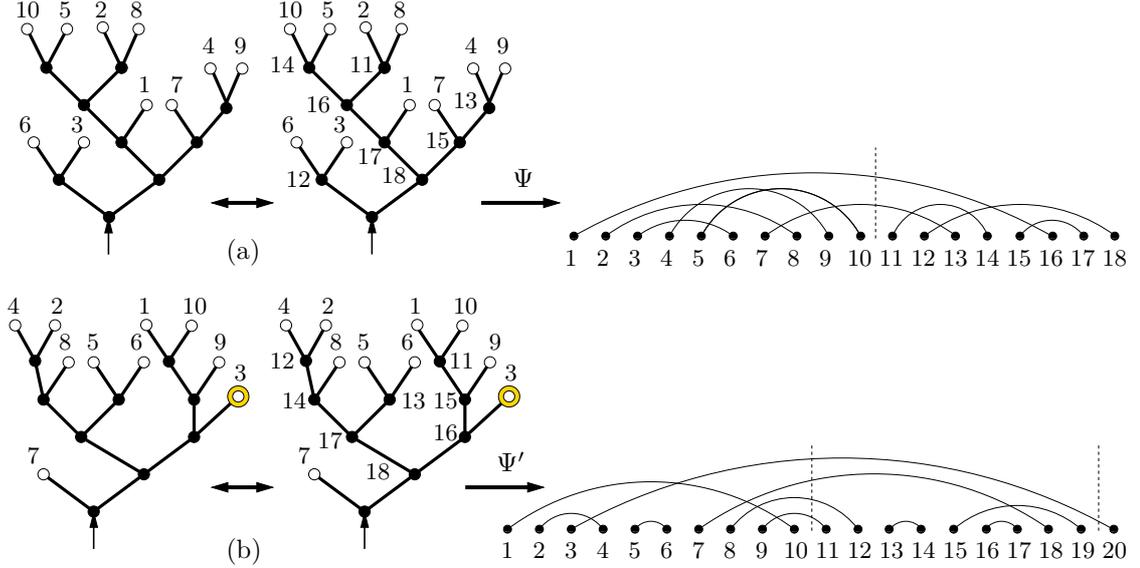}\caption{Bijections between labelled binary trees and matchings.}\label{fig:exp-matchings2} \end{center}\end{figure}

We will now present a bijection $\Psi$ due to Chen \cite{Chen:bij-algo-trees} between labelled binary trees of size $n$ and matchings on $[2n\!-\!2]$. We follow the exposition from \cite[p.16]{Stanley:volume2} for defining the bijection $\Psi$ and then define a similar bijection $\Psi'$ between marked binary trees of size $n$ and matchings on $[2n]$. The mappings $\Psi$ and $\Psi'$ are represented in Figure~\ref{fig:exp-matchings2}. The first step of this bijection is to attribute a label to each internal node of the binary tree.\\

\textbf{Definition of bijection $\Psi$ \cite{Chen:bij-algo-trees}}.
Let $\tau$ be a labelled binary tree with $n$ leaves labelled on $[n]$. One defines an induced
labelling of the $n\!-\!2$ internal non-root vertices of $\tau$ by the following procedure.
While there are unlabelled non-root vertices, we consider those among them that have both of their children labelled.
There is at least one such vertex and we consider the vertex $v$ which has the child with least label;
we then attribute to $v$ the least label in $[2n\!-\!2]\setminus[n]$ which has not yet been attributed. For instance, starting from the tree in Figure~\ref{fig:exp-labelled-trees}(a) one obtains the labels indicated in the \emph{fully} labelled tree represented in Figure~\ref{fig:exp-matchings2}(a). One then obtains the matching  $\pi=\Psi(\tau)$ on $[2n\!-\!2]$ by letting $\pi(i)=j$ for all pairs $i,j\in[2n-2]$ which are the labels of siblings in the fully labelled tree. In  \cite[p.16]{Stanley:volume2} the mapping $\Psi$ is proved to be a bijection between labelled binary trees of size $n$ and matchings on $[2n\!-\!2]$.\\

\textbf{Definition of bijection $\Psi'$}. Let $\tau'$ be a marked labelled binary tree of size $n$ and let $\tau$ be the corresponding unmarked tree. We consider the matching $\pi=\Psi(\tau)$ on $[2n\!-\!2]$ and define a matching $\pi'=\Psi'(\tau')$ on $[2n]$ as follows. If the marked vertex $v$ of $\tau'$ is the root, then $\pi'(i)=\pi(i)$ for all $i$ in $[2n\!-\!2]$ and $\pi'(2n-1)=2n$. If $v$ is not the root, then we consider its label $k$ and the label $l$ of its sibling. In this case $\pi'(i)=\pi(i)$ for all $i\neq k,l$ in $[2n\!-\!2]$, $\pi'(k)=2n$ and $\pi'(l)=2n-1$. It is clear that the mapping $\Psi'$ is a bijection between marked labelled binary tree of size $n$ and matchings on $[2n]$.\\

\smallskip

Combining the bijections $\Phi$, respectively $\Phi'$, with the bijections $\Psi$, respectively $\Psi'$, gives the following bijective proof of Proposition~\ref{prop:double-factorial}.

\begin{thm} The composition $\Psi\circ \Phi$ is a bijection between well-labelled Motzkin paths of size $n$ and matchings on $[2n\!-\!2]$. Similarly, $\Psi'\circ \Phi'$ is a bijection between well-labelled positive paths of size $n$ and matchings on~$[2n]$.\\
\end{thm}


\section{Enumerative corollaries} \label{section:counting}
We will now study the  number of horizontal steps in well-labelled paths through the bijections $\Phi$, $\Phi'$, $\Psi$, $\Psi'$ and deduce some enumerative corollaries in terms of the up-down sequences of permutations. Recall that a {\em horizontal step} of a well-labelled path $(\pp,\si)$ is a letter 0 in $\pp$.  We say that a leaf in a labelled binary tree is {\em single} if its sibling is an internal node.

\begin{thm}
 \label{thm:refbij1}
For all integers $n,k$, the mappings $\Phi$ and $\Psi$ induce successive bijections between \\
\ite well-labelled Motzkin paths of size $n$ with $k$ horizontal steps,\\
\ite labelled binary trees with $n$ leaves, $k$ of which are single leaves,\\
\ite matchings on $[2n-2]$ having $k$ pairs $(i,j)$ such that $i\in\{1,\ldots,n\}$ and $j\in\{n\!+\!1,\ldots,2n\!-\!2\}$.
\end{thm}

For example, the Motzkin path of size $n=10$ in Figure~\ref{fig:exp-labelled-path}(a) has 2 horizontal steps, the corresponding labelled binary tree represented in Figure \ref{fig:exp-labelled-trees}(a) has 2 single leaves, and the
corresponding matching represented in Figure~\ref{fig:exp-matchings2}(a) has 2 pairs $(i,j)$ such that  $i\in\{1,\ldots,n\}$ and $j\in\{n\!+\!1,\ldots,2n\!-\!2\}$.

\begin{proof}
\ite The correspondence between the number of horizontal steps of a Motzkin path $(\pp,\si)$ and the number of single leaves in the binary tree $\Phi(\pp,\si)$ follows from a simple induction on the size of $(\pp,\si)$. Indeed, one creates a single leaf in the recursive construction of $\Phi(\pp,\si)$ exactly when case (ii) thereof (corresponding to a horizontal step of $(\pp,\si)$) occurs.\\
\ite The correspondence between the number of single leaves in the binary tree $\tau$ and the number of pairs $(i,j)$ such that $i\in\{1,\ldots,n\}$ and $j\in\{n\!+\!1,\ldots,2n\!-\!2\}$ in the matching $\Psi(\tau)$ is an immediate consequence of the fact that the labels of external vertices are in $\{1,\ldots,n\}$ while the labels of internal vertices are in $\{n\!+\!1,\ldots,2n\!-\!2\}$.
 \end{proof}

\smallskip

\begin{cor}\label{cor:horizontal-Motzkin}
The number of well-labelled Motzkin paths of size $n$ having  $k$ horizontal steps is
$$a_{n,k}=  \binom{n}{k}\binom{n-2}{k}\,k!\,(n-k-1)!!\,(n-k-3)!!$$
if $n-k$ is even and 0 otherwise.
\end{cor}

\begin{proof}
By Theorem~\ref{thm:refbij1}, the number $a_{n,k}$ counts matchings on $[2n-2]$ with exactly $k$ pairs $(i,j)$ such that $i\in\{1,\ldots,n\}$ and $j\in\{n\!+\!1,\ldots,2n\!-\!2\}$. To enumerate such matchings, first choose these $k$ pairs: there are $\binom{n}{k}$ possibilities of choosing the integers $i$ in $\{1,\ldots,n\}$, there are $\binom{n-2}{k}$ possibilities for choosing the integers $j$ in $\{n\!+\!1,\ldots,2n\!-\!2\}$ and then $k!$ possibilities to define the pairing between the chosen integers in $\{1,\ldots,n\}$ and  the chosen integers in $\{n+1,\ldots,2n-2\}$.  After that, it remains to choose a pairing of the $n-k$ unmatched integers in $\{1,\ldots,n\}$ and a pairing of the $n-k-2$ unmatched integers in $\{n\!+\!1,\ldots,2n\!-\!2\}$. Such matchings exist only if $n-k$ is even and there are $(n-k-1)!!(n-k-3)!!$ of them in this case.
\end{proof}

\bigskip

We now examine horizontal steps in positive paths.  We say that a leaf in a marked labelled binary tree is {\em quasi-single} if it is not marked and its sibling is either marked or internal.

\begin{thm}
For all integers $n,k$, the mappings $\Phi'$ and $\Psi'$ induce successive bijections between \\
\ite well-labelled positive paths of size $n$ with $k$ horizontal steps,\\
\ite marked labelled binary trees of with $n$ leaves, $k$ of which are quasi-single leaves,\\
\ite matchings on $[2n]$ having $k$ pairs $(i,j)$ with $i\in \{1,\ldots,n\}$ and $j\in\{n+1,\ldots,2n-1\}$.
 \end{thm}

\begin{proof}
\ite  The correspondence between the number of horizontal steps of a positive path $(\pp,\si)$ and the number of quasi-single leaves in the marked tree $\Phi'(\pp,\si)$ follows from a simple induction on the size of $(\pp,\si)$.
Indeed, one creates a quasi-single leaf in the recursive construction of $\Phi'(\pp,\si)$ exactly when case (ii) of the definition of either $\Phi$ or $\Phi'$ occurs.\\
\ite We now consider a marked labelled binary tree $\tau'$ of size $n>1$. Let $v$ be the marked vertex and let $\tau$ be the non-marked tree obtained by forgetting the mark. The number $k$ of single leaves in $\tau$ and the number $k'$ of quasi-single leaves in $\tau'$ are related by
\begin{itemize}
\item $k'=k$ if $v$ is internal,
\item $k'=k-1$ if $v$ is a leaf and its sibling is internal
\item $k'=k+1$ if $v$ and its sibling are both leaves.
\end{itemize}
Similarly, the definition of $\Psi'$ gives a relation between the number $l'$ of pairs $(i,j)$ of the matching $\Psi'(\tau')$ such that $i\in\{1,\ldots,n\}$ and $j\in\{n+1,\ldots,2n-1\}$ and the number  $l$ of pairs $(i,j)$ of the matching $\Psi(\tau)$ such that $i\in\{1,\ldots,n\}$ and $j\in\{n+1,\ldots,2n-2\}$:
\begin{itemize}
\item $l'=l$ if the label of $v$ is larger than $n$,
\item $l'=l-1$ if the label of $v$ is not larger than $n$ and the label of its sibling is larger than~$n$,
\item $l'=l+1$ if the label of $v$ and its sibling are not larger than $n$.
\end{itemize}
Theorem~\ref{thm:refbij1} gives $k=l$, hence the previous discussion gives $k'=l'$ and concludes the proof.
\end{proof}

\smallskip

\begin{cor}\label{cor:horizontal-positive}
The number of well-labelled positive paths of size $n$ having $k$ horizontal steps is
\begin{equation}
b_{n,k}= \begin{cases}
\displaystyle~\binom{n}{k}\binom{n-1}{k}\, k!\,[(n-k-1)!!]^2& \textrm{if $n\!-\!k$ is even,}\\[10pt]
\displaystyle~ \binom{n}{k}\binom{n-1}{k}\,k!\,(n-k)!!\,(n-k-2)!!&\textrm{otherwise.}\\
  \end{cases}
\end{equation}\end{cor}

The proof of Corollary~\ref{cor:horizontal-positive} is very similar to the Corollary~\ref{cor:horizontal-Motzkin} and is omitted. We now study the consequence of these results in terms of the \emph{up-down sequences} of permutations. \\

An {\em ascent} of a permutation $\si=\si_1\si_2\ldots \si_n$ is an index $i<n$ such that $\si_i<\si_{i+1}$; a \emph{descent} is an index $i<n$ such that $\si_{i}>\si_{i+1}$. The enumeration of permutations with a given sequence of ascents and descents, called \emph{up-down sequences} (or \emph{shape}) was investigated for instance in \cite{Carlitz:prescribed-patterns,deBruijn:perm-given-updown,Foulkes:enum-perm-updown,Niven:comb-problem}. Here we will count permutations of size $n$ such that their up-down sequences belong to a certain family, while previous works focused on the enumeration according to a single up-down sequence. \\

We say that a permutation $\si$ has a \emph{positive up-down sequence} if for all $j\leq n$ the number of ascents less than $j$ is no more than the number of descents less than $j$.  We say that $\si$ has a \emph{Dyck up-down sequence} if it has one more ascent than descents but for all $j<n$ the number of ascents less than $j$ is no more than the number of descents less than $j$. Observe that a well-labelled path $(\pp,\si)$ with no horizontal steps is completely determined by the permutation $\si$ (indeed, the word $\pp$ is determined by the up-down sequence of $\si$). Moreover, the well-labelled path $(\pp,\si)$ is positive (resp. Motzkin) if and only if the permutation $\si$ has a positive (resp. Dyck) up-down sequence. Thus, the following theorem immediately follows by looking at the specialisation $k=0$ in Corollaries~\ref{cor:horizontal-Motzkin} and~\ref{cor:horizontal-positive}.

\begin{thm}
For any integer $n$, the number of permutations of size $n$ having a positive up-down sequence is $[(n-1)!!]^2$ if $n$ is even and $n!!\,(n-2)!!$ otherwise. The number of permutations of size $n$ having a Dyck up-down sequence is $(n-1)!!\, (n-3)!!$ if $n$ is even and 0 otherwise.
\end{thm}

We would be happy to see a more direct bijective proof of these specialisations.\\

\titre{Acknowledgement:} We are very thankful to Sylvie Corteel for fruitful discussions and for providing us with the conjectural formula for the numbers $b_{n,k}$.

\bibliography{biblio-Motzkinetiquetes}
\bibliographystyle{abbrv}

\end{document}